# CYCLES, SUBMANIFOLDS, AND STRUCTURES ON NORMAL BUNDLES

C. BOHR, B. HANKE, AND D. KOTSCHICK

ABSTRACT. We give explicit examples of degree 3 cohomology classes not Poincaré dual to submanifolds, and discuss the realisability of homology classes by submanifolds with Spin$^c$ normal bundles.

## 1. INTRODUCTION

In this paper we consider some classical questions concerning the realisation of cycles in smooth manifolds by submanifolds, in the spirit of Thom's celebrated paper [15]. We construct explicit examples to show that Thom's results are sharp in certain directions; in particular we give an explicit solution to a problem raised by Thom in [8]. We then discuss extensions to the case where one requires the normal bundle of the submanifold to have a Spin$^c$-structure.

Our interest in these questions was stimulated by two related developments in high-energy physics. One is the consideration of gerbes classified by elements of the third integral cohomology, the other is the appearance of the Spin$^c$ condition on normal bundles of certain submanifolds as an anomaly cancellation condition in string theory with D-branes.

Gerbes are topological objects classified by the third integral cohomology in the same way that complex line bundles are classified by the second integral cohomology. Chatterjee [3] and Hitchin [10] propose and pursue a down-to-earth approach to these objects, in which gerbes are constructed from codimension three submanifolds, rather than being thought of as sheaves of categories. The characteristic class of such a gerbe is the Poincaré dual of the homology class of the submanifold.

*Date*: October 30, 2018; MSC 2000: primary 57R95, secondary 81T30.

We are grateful to E. Scheidegger for explanations of the physical relevance of these problems. The second author would like to thank F. Cohen for useful discussions. The first two authors are supported by fellowships from the *Deutsche Forschungsgemeinschaft*; the second and third authors are members of the *European Differential Geometry Endeavour* (EDGE), Research Training Network HPRN-CT-2000-00101, supported by The European Human Potential Programme.

We will see in section 2 below that there are degree three cohomology classes which are not Poincaré dual to any submanifold, so that the construction via submanifolds does not yield all gerbes.

In [16] Witten postulated the existence of Spin$^c$-structures on the normal bundles of certain submanifolds arising in string theory. This condition was then interpreted as an anomaly cancellation condition [6]. Bryant and Sharpe [2] showed that in many of the physically relevant, low-dimensional cases, the assumption is satisfied automatically. In section 3 we exhibit examples which show that certain homology classes cannot be realised by submanifolds with Spin$^c$ normal bundle, although they can be realised by smooth submanifolds, and give a cohomological criterion for the existence of submanifolds with Spin$^c$ normal bundles realising a given homology class of degree 4.

In section 4 we show that there are homology classes with both kinds of representatives, with Spin$^c$ normal bundles and with non-Spin$^c$ normal bundles. These examples, or variations thereof, should be relevant to current work in M-theory, see [5].

Submanifolds with either almost complex or spin normal bundles have been discussed in the literature, see for example [13] for the almost complex case. Of course, almost complex or spin bundles are automatically Spin$^c$, but there are many homology classes which can be represented by submanifolds with Spin$^c$ normal bundles, although all their representatives have neither spin nor almost complex normal bundles. There are analogs of our results on Spin$^c$ normal bundles for the case of spin normal bundles. Generally speaking, the arguments to prove these are simplifications of what we do in the Spin$^c$ case. We leave it to the interested reader to work this out.

## 2. Realisation of cycles by submanifolds

We consider oriented smooth manifolds $M$ of dimension $n$, closed unless the contrary is stated. It is easy to see that every class in $H_1(M, \mathbb{Z})$, $H_{n-2}(M, \mathbb{Z})$ or $H_{n-1}(M, \mathbb{Z})$ can be represented by a smoothly embedded submanifold, and Hopf [11] proved that this is also the case for classes in $H_2(M, \mathbb{Z})$. This settles the representability of homology classes for manifolds of dimension at most 5.

Thom [15] introduced spaces $\mathrm{MSO}(k)$, which we now call Thom spaces, and which by construction have the following property: a class $x \in H_{n-k}(M, \mathbb{Z})$ can be realised by a smoothly embedded submanifold $\Sigma \subset M$ of codimension $k$ if and only if there is a map $f \colon M \to \mathrm{MSO}(k)$ which pulls back the Thom class $\tau \in H^k(\mathrm{MSO}(k), \mathbb{Z})$ to the Poincaré dual of $x$. Using this, he was able to show that all homology classes in degrees at most 6 are realised by submanifolds. The first case he left open was that of degree 7 classes in 10-dimensional manifolds, and this was included as Problem 8 in [8]. Thom did give an example of a manifold of dimension 14 with a homology class of degree 7 that cannot be realised by a submanifold, see [15], page 62.

**Theorem 1.** *For every $n \geq 10$ there are $n$-dimensional manifolds $M$ and homology classes $x \in H_{n-3}(M, \mathbb{Z})$ that cannot be realised by submanifolds. The manifolds $M$ can be chosen to have torsion-free cohomology.*

*Proof.* Consider the symplectic group

$$\mathrm{Sp}(2) = \{ A \in \mathbb{H}^2 \mid AA^* = 1 \} \,,$$

where $A^*$ denotes the quaternionic conjugate of the transpose of $A$. This is a simply connected simple Lie group of dimension 10, and the standard inclusion $\mathrm{Sp}(1) \subset \mathrm{Sp}(2)$ induces a fibre bundle

$$S^3 \hookrightarrow \mathrm{Sp}(2) \to S^7 \,.$$

The integral cohomology of $\mathrm{Sp}(2)$ is isomorphic to $\mathbb{Z}$ in dimensions 0, 3, 7 and 10, and is zero in all other dimensions. It is well-known, see for example [1], that $\mathcal{P}_3^1$, the first reduced Steenrod power operation for the prime number 3, induces an isomorphism

$$H^3(\mathrm{Sp}(2); \mathbb{Z}_3) \cong H^7(\mathrm{Sp}(2); \mathbb{Z}_3) \,.$$

By Poincaré duality, the cup product pairing

$$H^3(\mathrm{Sp}(2); \mathbb{Z}_3) \times H^7(\mathrm{Sp}(2); \mathbb{Z}_3) \to H^{10}(\mathrm{Sp}(2); \mathbb{Z}_3)$$

is nontrivial. Let $x \in H_7(\mathrm{Sp}(2); \mathbb{Z})$ be a generator, $u \in H^3(\mathrm{Sp}(2); \mathbb{Z})$ its Poincaré dual (with respect to some orientation of $\mathrm{Sp}(2)$) and suppose there is a map

$$f \colon \mathrm{Sp}(2) \longrightarrow \mathrm{MSO}(3)$$

with $f^*(\tau) = u$, where $\tau \in H^3(\mathrm{MSO}(3); \mathbb{Z})$ is the Thom class. If we denote by $\overline{\tau} \in H^3(\mathrm{MSO}(3); \mathbb{Z}_3)$ the reduction of $\tau$, then, by the above considerations, we obtain

$$f^*(\overline{\tau} \cup \mathcal{P}_3^1(\overline{\tau})) = f^*(\overline{\tau}) \cup \mathcal{P}_3^1(f^*(\overline{\tau})) \neq 0 \in H^{10}(\mathrm{Sp}(2); \mathbb{Z}_3) \,.$$

But $H^*(\mathrm{MSO}(3); \mathbb{Z}_3)$ is concentrated in degrees $3 + 4k$ by the Thom isomorphism and the well-known fact that $H^*(\mathrm{BSO}(3); \mathbb{Z}_3)$ is a polynomial ring on the mod 3 reduction of the universal first Pontryagin class. This is a contradiction. Hence there is no $f$ as above, and $x$ cannot be realised by a submanifold.

Higher dimensional examples can be constructed by taking products of $\mathrm{Sp}(2)$ with spheres, say.                                       $\square$

*Remark* 1. As $\pi_6(S^3) \cong \mathbb{Z}_{12}$, there are 12 different principal $S^3$-bundles over $S^7$. All the total spaces have the same cohomology rings, and for the bundles corresponding to elements in $\mathbb{Z}_{12}$ whose orders are not powers of 2 one can prove that the generator of the degree 7 homology is not realisable by a submanifold using $\mathcal{P}_3^1$ as in the above proof.

*Remark* 2. In the first version of this paper we constructed homology classes of codimension 3 not realisable by submanifolds in manifolds of dimension $\geq 10$ which are products of lens spaces. That construction is an improvement of the construction of Thom [15], page 62, and relies on the presence of torsion in integral homology. The above examples are more interesting because their integral homology is torsion-free. In general, it is well-known (cf. [4]) that for an arbitrary manifold $M$ with torsion-free integral homology the Atiyah-Hirzebruch spectral sequence converging to the oriented bordism group $\Omega_*^{\mathrm{SO}}(M)$ collapses at the $E_2$-level. This implies that all homology classes in $H_*(M; \mathbb{Z})$ can in fact be realised by immersed submanifolds. Up to the middle dimension one can then choose the submanifolds to be embedded.

*Remark* 3. Problem 8 in [8] was also solved by Kreck, see [9], page 771. His argument is indirect, using the Atiyah–Hirzebruch spectral sequence and does not yield any explicit examples.

Instead of just asking for some oriented submanifold representing a given homology class, it is natural to look for connected submanifolds. It is easy to see that in codimension one a non-zero class is representable by a connected submanifold if and only if it is primitive, or indivisible, in integral homology, and that in codimension $\geq 2$ every class that can be represented by a submanifold can also be represented by one that is connected.

In the rest of this paper we shall discuss representations by submanifolds with $\mathrm{Spin}^c$ normal bundles. In codimension one this is not an issue, because an oriented codimension one submanifold of an oriented manifold has trivial normal bundle. In higher codimension we can always assume the representatives to be connected. We then have the following obvious lemma.

**Lemma 1.** *Let $M$ be an oriented manifold, and $\Sigma \subset M$ an oriented submanifold. Assume that the restriction $TM|_\Sigma$ is $\mathrm{Spin}^c$. Then $\Sigma$ is $\mathrm{Spin}^c$ if and only if its normal bundle is.*

The assumption is satisfied in particular if $M$ is $\mathrm{Spin}^c$ itself.

At various points we shall use Whitney's classical result that all smooth oriented 4-manifolds are $\mathrm{Spin}^c$.

### 3. Homological obstructions to Spin$^c$-structures on normal bundles

We now discuss the representation of homology classes by submanifolds with Spin$^c$ normal bundles. As an oriented vector bundle $E$ over an oriented manifold is Spin$^c$ as soon as $\mathrm{rk}(E) \leq 2$ or $\dim(M) \leq 3$, the first interesting case is that of degree 4 homology classes in 7-manifolds. By the results of Thom [15], all these classes are representable by submanifolds, but we shall see that the Spin$^c$ condition on the normal bundle is a non-trivial constraint.

Our first main result in this section, Theorem 2 below, gives a necessary and sufficient cohomological condition for the existence of a 4-dimensional oriented submanifold $X$ of a manifold $M$ with Spin$^c$ normal bundle representing a given class $\alpha \in H_4(M, \mathbb{Z})$. To formulate this criterion, we first have to recall the definition of the integral third Steenrod square. Changing notation from the previous section, let

$$\overline{\beta} \colon H^*(M, \mathbb{Z}_2) \longrightarrow H^{*+1}(M, \mathbb{Z})$$

denote the Bockstein operator associated to the exact sequence

$$0 \to \mathbb{Z} \to \mathbb{Z} \to \mathbb{Z}_2 \to 0 \ .$$

Then the integral Steenrod square $\overline{Sq^3}$ is defined to be the cohomology operation given by the composition $\overline{\beta} \circ Sq^2 \circ \mu_2$:

$$\overline{Sq^3} \ : \ H^*(M, \mathbb{Z}) \xrightarrow{\mu_2} H^*(M, \mathbb{Z}_2) \xrightarrow{Sq^2} H^{*+2}(M, \mathbb{Z}_2) \xrightarrow{\overline{\beta}} H^{*+3}(M, \mathbb{Z}) \ .$$

Here, $\mu_2$ denotes reduction mod 2. Note that the image of this operation is contained in the 2-torsion, and that its reduction mod 2 is the reduction followed by the usual Steenrod square $Sq^3 = Sq^1 \circ Sq^2$.

**Lemma 2.** *Let $M$ be an oriented closed manifold of dimension $n$ and $i\colon X \hookrightarrow M$ a connected oriented $m$-dimensional submanifold. Let $\pi\colon N \to X$ denote the normal bundle of $X$ in $M$ (considered as a tubular neighbourhood of $X$) and*

$$f\colon H^*(X, \mathbb{Z}) \cong H^{*+(n-m)}(N, \partial N, \mathbb{Z}) \cong H^{*+(n-m)}(M, M \setminus N, \mathbb{Z})$$
$$\longrightarrow H^{*+(n-m)}(M, \mathbb{Z})$$

*the composition of the Thom isomorphism with the excision map. If $PD_X$ and $PD_M$ denote the Poincaré duality isomorphisms associated*

*to $X$ and $M$ respectively, then the following diagram is commutative:*

$$
\begin{array}{ccc}
H^*(X,\mathbb{Z}) & \xrightarrow{\ f\ } & H^{*+(n-m)}(M,\mathbb{Z}) \\
PD_X \downarrow & & PD_M \downarrow \\
H_{m-*}(X,\mathbb{Z}) & \xrightarrow{\ i_*\ } & H_{m-*}(M,\mathbb{Z})\,.
\end{array}
$$

*Proof.* Let $\tau \in H^{n-m}(N,\partial N,\mathbb{Z})$ denote the Thom class (induced by the orientations of $X$ and $M$) and let $[X] \in H_m(X,\mathbb{Z})$, $[N,\partial N] \in H_n(N,\partial N,\mathbb{Z})$ and $[M] \in H_n(M,\mathbb{Z})$ be the fundamental classes. We shall denote by $j_*\colon H_*(X,\mathbb{Z}) \to H_*(N,\mathbb{Z})$ and $k_*\colon H_*(N,\mathbb{Z}) \to H_*(M,\mathbb{Z})$ the maps induced by inclusion. Note that $\tau \cap [N,\partial N] = j_*([X])$. For $\varphi \in H^*(X,\mathbb{Z})$ the claim is shown by the following calculation in $H_{m-*}(M,\mathbb{Z})$:

$$
\begin{aligned}
f(\phi) \cap [M] &= k_*((\pi^*(\phi) \cup \tau) \cap [N,\partial N]) \\
&= k_*(\pi^*(\phi) \cap (\tau \cap [N,\partial N])) \\
&= k_*(\pi^*(\phi) \cap j_*([X])) \\
&= k_*(j_*(j^*(\pi^*(\phi)) \cap [X])) \\
&= k_*(j_*(PD_X(\phi))) \\
&= i_*(PD_X(\phi))\,.
\end{aligned}
$$

$\square$

We denote by $W_3$ the third integral Stiefel-Whitney class, i. e. the image of $w_2$ under the Bockstein map $\overline{\beta}$.

**Lemma 3.** *Let $M$ be an oriented closed manifold of dimension $n \geq 5$ and $i\colon X \hookrightarrow M$ the inclusion of a connected oriented 4-dimensional submanifold. Denote $\varphi = PD_M(i_*[X])$. Then*

$$
i_*(PD_X(W_3(M)|_X)) = PD_M(\overline{Sq^3}(\varphi)),
$$

*where $\overline{Sq^3}$ is the integral third Steenrod square defined above.*

*Proof.* Let $\nu$ denote the normal bundle of $X$ in $M$. Since the 4-manifold $X$ is Spin$^c$, we have $W_3(M)|_X = W_3(\nu)$. According to a result of Thom [14], the Thom isomorphism maps $W_3(\nu)$ to $\overline{Sq^3}(\tau)$, where $\tau$ denotes the Thom class. Using the notations from Lemma 2, this implies that $f(W_3(M)|_X) = \overline{Sq^3}(\varphi)$, and the claim follows from Lemma 2. $\square$

We are now ready to prove the following result.

**Theorem 2.** *Let $M$ be a smooth closed oriented manifold of dimension $n \geq 7$, and $\alpha \in H_4(M,\mathbb{Z})$ a homology class with Poincaré dual $\varphi \in H^{n-4}(M,\mathbb{Z})$. Then $\alpha$ can be realised by an embedded submanifold with* Spin$^c$ *normal bundle if and only if $\overline{Sq^3}(\varphi) = 0$.*

*Proof.* It is an immediate consequence of Lemma 3 that the condition is necessary, because for an embedded 4-manifold $X$ with normal bundle $\nu$, we have $W_3(M)|_X = W_3(\nu)$.

Conversely, suppose that $\overline{Sq^3}(\varphi) = 0$. By Theorem II.27 in [15], there exists an embedded 4-manifold $X \subset M$ representing the homology class $\alpha$. Let $S \subset X$ be a circle dual to the restriction $W_3(M)|_X$. By Lemma 3, the homology class of $S$ in $M$ is dual to $\overline{Sq^3}(\varphi) = 0$, hence $S$ bounds a surface $F \subset M$. We can arrange that $F$ meets $X$ only in $\partial F = S$. Now we will apply ambient surgery along $S$, i. e. we will replace a tubular neighborhood of $S$ in $X$ by $Y = F \times S^2$, as follows.

The normal bundle of $F$ in $M$ is trivial. If we pick a trivialisation, then the normal bundle of $S$ in $X$ is a rank 3 subbundle of this trivialised bundle, i. e. it is given by a map of $S$ into the Grassmannian $G_3(\mathbb{R}^{n-2})$ of oriented 3-dimensional subspaces of $\mathbb{R}^{n-2}$. A framing of the circle $S$ corresponds to a lift of this map to the Stiefel manifold $V_3(\mathbb{R}^{n-2})$ of 3-frames in $\mathbb{R}^{n-2}$, and since $V_3(\mathbb{R}^{n-2})$ is simply connected, we can extend the lift to a map $f \colon F \to V_3(\mathbb{R}^{n-2})$. This map in turn defines a trivial rank 3 subbundle of the normal bundle of $F$, and if we use this bundle, respectively its trivialisation given by $f$, to thicken up $F$, we obtain an embedding of $F \times D^3$ in $M$ such that $(F \times D^3) \cap X = S \times D^3$ is a tubular neighbourhood $T$ of $S$. If we let $Y = F \times S^2 \subset F \times D^3$ and $X_0 = X \setminus T$, then we can replace $T$ by $Y$, i. e. we can form the manifold $X' = X_0 \cup_{\partial Y} Y \subset M$, which still represents the class $\alpha$.

We will compute the restriction of $W_3(M)$ to $X'$. For the remainder of this proof, all homology and cohomology groups are understood to have integral coefficients.

The commutative diagram

$$
\begin{array}{ccc}
H^3(X, X_0) & \longrightarrow & H^3(X) \\
\downarrow & & \downarrow \\
H_1(T) & \longrightarrow & H_1(X)
\end{array}
$$

in which the vertical arrows are given by duality in $X$, shows that the image of $H^3(X, X_0) = \mathbb{Z}$ in $H^3(X)$ is generated by $W_3(M)|_X$. Hence, by exactness of the sequence

$$ H^3(X, X_0) \longrightarrow H^3(X) \longrightarrow H^3(X_0), $$

we get $W_3(M)|_{X_0} = 0$.

By construction, $X'$ decomposes as $X' = X_0 \cup_{\partial Y} Y$, hence we can use the Mayer–Vietoris sequence

$$ H^2(X_0) \oplus H^2(Y) \longrightarrow H^2(\partial Y) \xrightarrow{\ \partial\ } H^3(X') \to H^3(X_0) \oplus H^3(Y) \ . $$

The 2-torsion class $W_3(M)$ restricts to zero on $Y$, since $H^3(Y)$ is torsion–free. Using this together with $W_3(M)|_{X_0} = 0$, we obtain $W_3(M)|_{X'} = \partial(x)$ for some $x \in H^2(\partial Y)$. However, since the restriction $H^2(Y) \to H^2(\partial Y)$ is onto, so is $H^2(X_0) \oplus H^2(Y) \to H^2(\partial Y)$, and hence $\partial = 0$. This proves that $W_3(M)|_{X'} = 0$. By Lemma 1 and using that every orientable 4-manifold admits a Spin$^c$-structure, we conclude that the normal bundle of $X'$ is Spin$^c$, as desired. $\qquad\square$

*Remark* 4. For $n \geq 9$ we will sketch a homotopy-theoretic proof in the style of Thom [15] for the sufficiency of the cohomological condition in Theorem 2.

Consider a map

$$f \colon \mathrm{MSpin}^c(n-4) \to K(\mathbb{Z}, n-4)$$

pulling back the canonical class in $H^{n-4}(K(\mathbb{Z}, n-4), \mathbb{Z})$ to the Thom class $\tau \in H^{n-4}(\mathrm{MSpin}^c(n-4), \mathbb{Z})$. It is easy to see that this is an $(n-2)$-equivalence whose homotopy fibre $F$ has $\pi_{n-2}(F) \cong \mathbb{Z}$. Hence, by standard arguments in obstruction theory, there exists a fibration $K \to K(\mathbb{Z}, n-4)$ with fiber $K(\mathbb{Z}, n-2)$ and a lift

$$g \colon \mathrm{MSpin}^c(n-4) \to K$$

of $f$, which is an $(n-1)$-equivalence. We can therefore construct a map $g' \colon K^{(n-1)} \to \mathrm{MSpin}^c(n-4)$ defined on the $(n-1)$-skeleton of $K$, which, in cohomology, pulls back the Thom class in $H^{n-4}(\mathrm{MSpin}^c(n-4), \mathbb{Z})$ to the canonical class in $H^{n-4}(K(\mathbb{Z}, n-4), \mathbb{Z}) = H^{n-4}(K, \mathbb{Z})$. Now, suppose that $\varphi \in H^{n-4}(M, \mathbb{Z})$ and $\overline{Sq^3}(\varphi) = 0$. The map $M \to K(\mathbb{Z}, n-4)$ corresponding to $\varphi$ can be lifted to a map $h \colon M \to K$, because the obstruction for such a lift is $\overline{Sq^3}(\varphi)$, which vanishes by assumption. Using the previous argument and cellular approximation, one gets a map

$$M^{(n-1)} \to \mathrm{MSpin}^c(n-4)$$

which, in cohomology, pulls back the Thom class of $\mathrm{MSpin}^c(n-4)$ to $\varphi$. Extension of this map to all of $M$ is possible, if an obstruction class in $H^n(M, \pi_{n-1}(\mathrm{MSpin}^c(n-4)))$ vanishes. But in the stable range $n \geq 9$ we have

$$\pi_{n-1}(\mathrm{MSpin}^c(n-4)) = \Omega_3^{\mathrm{Spin}^c} = 0$$

by the results in [12].

We shall now describe two explicit constructions yielding examples of homology classes which can be represented only by submanifolds with non-Spin$^c$ normal bundles. One of these constructions, given in Theorem 3, yields degree 4 classes in $n$-dimensional manifolds. The other construction, in Theorem 4, gives codimension 3 classes in $n$–dimensional manifolds. In both cases, $n$ can take all values $\geq 7$. For $n = 7$ the two constructions are essentially the same.

**Theorem 3.** *For every $n \geq 7$ there exists an $n$–dimensional manifold $M$ and a homology class $\alpha \in H_4(M, \mathbb{Z})$ for which every embedded representative $Y \subset M$ has non-*Spin$^c$ *normal bundle.*

*In fact, the conclusion holds for all embedded $Y$ with $[Y] \equiv \alpha \mod 2$.*

Note that by the results of Thom [15] all degree 4 classes are representable by submanifolds.

**Theorem 4.** *For every $n \geq 7$ there exists an $n$–dimensional manifold $M$ and a homology class $\alpha \in H_{n-3}(M, \mathbb{Z})$ representable by an embedded submanifold, for which every such representative has non-*Spin$^c$ *normal bundle.*

*In fact, the conclusion holds for all embedded $Y$ with $[Y] \equiv \alpha \mod 2$.*

To prove these results, we need three more Lemmata.

**Lemma 4.** *Let $X$ be a 4–manifold and $x \in H^2(X, \mathbb{Z}_2)$ an arbitrary class. Then for every $n \geq 3$ there exists an orientable rank $n$ vector bundle $E \to X$ with $w_2(E) = x$.*

*Proof.* Since $\pi_1(\mathrm{BSO}(3)) = 0 = \pi_3(\mathrm{BSO}(3))$ and $\pi_2(\mathrm{BSO}(3))) = \mathbb{Z}_2$, the map $h \colon \mathrm{BSO}(3) \to K(\mathbb{Z}_2, 2)$ representing $w_2 \in H^2(\mathrm{BSO}(3), \mathbb{Z}_2)$ is a 4-equivalence. Hence, there is a map $f \colon X \to \mathrm{BSO}(3)$, such that a map $X \to K(\mathbb{Z}_2, 2)$ realising $x$ can be written as $h \circ f$ (using cellular approximation). Therefore $f^* w_2 = x$.

If $n = 3$, we take the pullback of the universal bundle by $f$. For larger $n$ we can add a trivial bundle. □

**Lemma 5.** *Let $X$ be a 4–manifold with $H_1(X, \mathbb{Z}) = \mathbb{Z}_2$, and $\pi \colon M \to X$ a $k$–sphere bundle with a section $\sigma \colon X \to M$. Assume $k \geq 1$. Let $j \colon Y \hookrightarrow M$ be an embedding of another oriented connected 4–manifold $Y$ with $\sigma_*[X] \equiv j_*[Y] \mod 2$. Denote by $f \colon Y \to X$ the composition $\pi \circ j$.*

*Then $f^* \colon H^3(X, \mathbb{Z}) \to H^3(Y, \mathbb{Z})$ is injective.*

*Proof.* Using Poincaré duality and the Bockstein exact sequence, one easily checks that $H^3(X, \mathbb{Z}) \to H^3(X, \mathbb{Z}_2)$ is an isomorphism. Hence it suffices to show that $f^* \colon H^3(X, \mathbb{Z}_2) \to H^3(Y, \mathbb{Z}_2)$ is injective. Therefore we work with $\mathbb{Z}_2$ coefficients throughout the proof.

The mod 2 mapping degree $\deg_2(f)$ of $f$ is nothing but the mod 2 intersection number $S \cdot j(Y)$ of a fiber $S$ with $j(Y)$. Since this number can be expressed as the cup product of the classes dual to $j_*[Y]$ and $[S]$, it only depends on the $\mathbb{Z}_2$–homology class of $Y$, and hence

$$\deg_2(f) = S \cdot j(Y) = S \cdot \sigma(X) = 1 \in \mathbb{Z}_2 \ .$$

Hence $f^*[X] = [Y]$ in cohomology with $\mathbb{Z}_2$ coefficients.

Now assume $x \in H^3(X, \mathbb{Z}_2)$ is not zero. By duality, there exists a class $y \in H^1(X, \mathbb{Z}_2)$ such that $x \cup y = [X] \in H^4(X, \mathbb{Z}_2) = \mathbb{Z}_2$. We then have

$$f^*x \cup f^*y = f^*(x \cup y) = f^*[X] = [Y] \neq 0 \ ,$$

and this implies $f^*x \neq 0$, proving that $H^3(X, \mathbb{Z}_2) \to H^3(Y, \mathbb{Z}_2)$ is indeed injective. $\qquad\square$

**Lemma 6.** *Let $M$ be a manifold and $\varphi \in H^3(M, \mathbb{Z})$ a cohomology class whose Poincaré dual can be represented by a submanifold with Spin$^c$ normal bundle. Then $\varphi \cup \varphi = 0$.*

*Proof.* If $\tau \in H^3(\mathrm{MSpin}^c(3), \mathbb{Z})$ denotes the Thom class, $W_3(\gamma)$ is mapped to $\tau \cup \tau$ under the Thom isomorphism, where $\gamma$ denotes the canonical bundle over $BSpin^c(3)$. As the bundle $\gamma$ is Spin$^c$, we get $W_3(\gamma) = 0$ and therefore $\tau \cup \tau = 0$. Let $f \colon M \to \mathrm{MSpin}^c(3)$ be a map with $f^*(\tau) = \varphi$. This implies $\varphi \cup \varphi = f^*(\tau \cup \tau) = 0$ as claimed. $\qquad\square$

*Proof of Theorem 3.* Choose a 4–manifold $X$ with $H_1(X, \mathbb{Z}) = \mathbb{Z}_2$. For convenience we can take $X$ to be spin.

We have $H^3(X, \mathbb{Z}) = \mathbb{Z}_2$ and $\beta \colon H^2(X, \mathbb{Z}_2) \to H^3(X, \mathbb{Z})$ is surjective. Choose a class $x \in H^2(X, \mathbb{Z}_2)$ with $\beta(x) \neq 0$. By Lemma 4, there exists a vector bundle $E \to X$ of rank $n - 4$ with $w_2(E) = x$, and hence $W_3(E) = \beta(w_2(E)) \neq 0$. So $E \to X$ is not Spin$^c$.

Now let $M$ denote the double of the disk bundle of $E$. Then $M$ is an $S^{n-4}$–bundle over $X$, in fact it is the sphere bundle of $E \oplus \mathbb{R}$. Let $\pi \colon M \to X$ denote the projection and $\sigma \colon X \to M$ the section given by the zero section of $E$. Let $\alpha = \sigma_*[X] \in H_4(M, \mathbb{Z})$. By construction, $\alpha$ is represented by the submanifold $X' = \sigma(X)$. We could now use Theorem 2 to show that the class $\alpha$ can not be represented by a submanifold having Spin$^c$ normal bundle. Alternatively, we can argue as follows.

We claim that $w_2(M) = \pi^*w_2(E)$. In fact, the Gysin sequence shows that the restriction $H^2(M, \mathbb{Z}_2) \to H^2(X', \mathbb{Z}_2)$ is an isomorphism, so we only have to show that $w_2(M)$ and $\pi^*w_2(E)$ have the same restriction to $X'$. Let $\nu(X')$ denote the normal bundle of $X'$ in $M$. This is just $(\pi|_{X'})^*E$. Now

$$w_2(M)|_{X'} = w_2(\nu(X')) = (\pi|_{X'})^*w_2(E) = \pi^*w_2(E)|_{X'} \ ,$$

where the first equality holds since $TM|_{X'} = TX' \oplus \nu(X')$ and $X'$ is spin. This proves our claim.

Applying the Bockstein, we also obtain $W_3(M) = \pi^* W_3(E)$. Now suppose that $j\colon Y \to M$ is an embedding of a 4–manifold $Y$ with $j_*[Y] \equiv \alpha \mod 2$. Consider the composition $f = \pi \circ j\colon Y \to X$. By Lemma 5, the induced map $f^*\colon H^3(X, \mathbb{Z}) \to H^3(Y, \mathbb{Z})$ is injective. But, since $Y$ is Spin$^c$, we have

$$W_3(\nu(Y)) = W_3(M)|_Y = j^*\pi^* W_3(E) = f^* W_3(E) \neq 0,$$

and hence the normal bundle $\nu(Y)$ is not Spin$^c$. $\qquad\qquad\square$

*Proof of Theorem 4.* Choose a connected oriented 4–manifold $Y$ whose first homology $H_1(Y, \mathbb{Z})$ has non–trivial 2-torsion. From the Bockstein exact sequence it is clear that there exists a class in $H^2(Y, \mathbb{Z}_2)$ which cannot be lifted to $H^2(Y, \mathbb{Z})$. By Lemma 4, this class can be realised as the second Stiefel–Whitney class of an orientable rank 3 vector bundle $L \to Y$, which is therefore not Spin$^c$. Pick a connected $(n-7)$–manifold $Y'$ and let $X = Y \times Y'$. Let $E \to X$ denote the pullback of $L$ to $X$. Then $E$ is not Spin$^c$, since its restriction to a factor is not Spin$^c$.

Denote by $M$ the double of the disk bundle of $E$. Then $\pi\colon M \to X$ is an $S^3$–bundle. As in the proof of Theorem 3, we have a section $\sigma\colon X \to M$ given by the zero section in $E$, and the image $\sigma(X) \subset M$ has normal bundle $E$.

Let $\alpha = \sigma_*[X] \in H_{n-3}(M, \mathbb{Z})$ and $\varphi = PD(\alpha) \in H^3(M, \mathbb{Z})$. By construction, the class $\alpha$ is representable by an embedded submanifold. By Lemma 6, it is now sufficient to show that $\varphi' \cup \varphi' \neq 0$ for every class $\varphi' \in H^3(M, \mathbb{Z})$ with $\varphi' \equiv \varphi \mod 2$. For every such class, there exists an $x \in H^3(M, \mathbb{Z})$ with $\varphi' = \varphi + 2x$. Since we are in odd degree, $\varphi' \cup \varphi' = \varphi \cup \varphi$. Hence we only have to show that $\varphi \cup \varphi$ does not vanish.

Let $X' \subset M$ denote a copy of $X$ embedded by a section transverse to the zero section in the normal bundle of $\sigma(X)$. The intersection $\sigma(X) \cap X'$ is then a submanifold of dimension $n-6$ whose homology class is Poincaré dual to $\varphi \cup \varphi$. Since the homomorphism $H_*(\sigma(X), \mathbb{Z}) \to H_*(M, \mathbb{Z})$ induced by the inclusion is injective (a left inverse is given by $\pi_*$), the Theorem is proved if we can show that the homology class of $\sigma(X) \cap X'$ in $\sigma(X)$ is not zero. But this class in $H_{n-6}(\sigma(X), \mathbb{Z})$ is dual to the Euler class of the normal bundle of $\sigma(X)$, so if we identify $\sigma(X)$ with $X$, it is dual to the Euler class $e(E) \in H^3(X, \mathbb{Z})$. However, for an orientable rank 3 bundle, the Euler class is just the integral Stiefel–Whitney class $W_3$, and hence $e(E) = W_3(E) \neq 0$, since $E$ is not Spin$^c$. $\qquad\qquad\square$

## 4. Further examples

In Section 3 we saw examples of homology classes which do not have any representative submanifolds with Spin$^c$ normal bundle. On the other hand, there are plenty of examples where each representative must have Spin$^c$ normal bundle, for dimension reasons, say. In this section we show that most homology classes do not fall into one of these extreme cases, in that they have both kinds of representative submanifolds, some with Spin$^c$ normal bundle, and some with non-Spin$^c$ normal bundle. Put differently, the existence of a Spin$^c$-structure on the normal bundle of a submanifold is not usually a homological property.

The next example shows that sometimes a representative with non-Spin$^c$ normal bundle can be modified to one with Spin$^c$ normal bundle.

*Example* 1. We start as in the proof of Theorem 3. Suppose that $X$ is a 4–dimensional manifold with 2–torsion in $H_1(X, \mathbb{Z})$. Pick an oriented vector bundle $E \rightarrow X$ of rank 3 which does not admit a Spin$^c$–structure and let $M$ denote the double of the disk bundle of $E$. We then have an embedding $X \hookrightarrow M$ and the normal bundle of $X$ in $M$ is precisely $E$, in particular it is not Spin$^c$.

Now we can find a finite collection of pairwise disjoint circles in $M$, disjoint from $X$, such that their homotopy classes normally generate $\pi_1(M)$. Performing surgery along these circles yields a simply connected 7–manifold $N$ which still contains $X$ as a submanifold. Since we did not change anything in a neighborhood of $X$, the normal bundle of $X$ in $N$ is still $E$. Let $\alpha \in H_4(N, \mathbb{Z})$ denote the homology class of $X$. Then $\overline{Sq^3}(PD(\alpha)) = 0$, simply because the group $H^6(N, \mathbb{Z}) = H_1(N, \mathbb{Z}) = 0$ by construction, and hence Theorem 2 shows that there exists a second representative $X'$ of $\alpha$ whose normal bundle is Spin$^c$. Also note that this representative can be obtained from $X$ by ambient surgery as in the proof of Theorem 2. Since the manifold $N$ is simply connected, we could also use ambient surgery to produce a simply connected representative on which every orientable vector bundle is Spin$^c$, since its homology is torsion–free.

In the other direction, a representative with Spin$^c$ normal bundle can often be corrupted, to obtain one with non-Spin$^c$ normal bundle:

**Theorem 5.** *Let $M$ be a closed oriented $n$-manifold and $X \subset M$ an embedded $i$-manifold with* Spin$^c$ *normal bundle. If $n \geq 2i - 1 \geq 9$, then there is a submanifold $X' \subset M$ homologous to $X$ with non-*Spin$^c$ *normal bundle.*

In the proof we shall make use of the following:

**Lemma 7.** *For every $i \geq 5$ there are closed oriented $i$-manifolds which are not* Spin$^c$.

*Proof.* It suffices to prove the case $i = 5$. A 5-manifold $X$ which is not Spin$^c$ was exhibited in [6]. It is the total space of a $\mathbb{C}P^2$-bundle over $S^1$ whose monodromy is complex conjugation.

Another example is the symmetric space $Z = \mathrm{SU}(3)/\mathrm{SO}(3)$, where the inclusion is the canonical one. As a smooth manifold, this can be obtained from the $\mathbb{C}P^2$-bundle $X$ above by surgery killing the generator of $\pi_1(X)$. Clearly this does not affect the property of not being Spin$^c$. □

*Proof of Theorem 5.* By the above Lemma, there is a closed oriented $i$-manifold $Y$ which is not Spin$^c$. As we assumed $n \geq 2i - 1$, $Y$ can be embedded in $S^n$ by the results of [7]. By Lemma 1 its normal bundle is not Spin$^c$.

Now form the internal connected sum $X' = X \# Y$ in $M = M \# S^n$. It is clear that $X'$ is homologous to $X$. We have an inclusion $Y \setminus B \hookrightarrow X'$, where $B$ is a ball in $Y$, and since the inclusion $Y \setminus B \hookrightarrow Y$ induces an isomorphism in $H^3$, the restriction of the normal bundle of $X'$ to $Y \setminus B$ – which coincides with the restriction of the normal bundle of $Y$ in $S^n$ to $Y \setminus B$ – has non–vanishing $W_3$. Hence the normal bundle of $X'$ is not Spin$^c$, as desired. □

*Example* 2. Rees and Thomas [13] have given examples of homology classes realised by singular algebraic cycles in smooth algebraic varieties which cannot be realised by smooth submanifolds with almost complex normal bundles. The simplest example is the class of $[T_1] + [T_2]$ in $T_1 \times T_2$, where the $T_i$ are Abelian varieties of complex dimension 2. This is clearly representable by smoothly embedded 4-manifolds $X$. As every such $X$ is Spin$^c$, Lemma 1 shows that the normal bundle is also Spin$^c$, although it cannot be almost complex.

Returning now to the motivating situation in physics, Witten [16] considered Spin$^c$-structures on the normal bundles of Lagrangian submanifolds in Kähler manifolds, particularly Calabi-Yaus. In this case the ambient manifold is Spin$^c$ and Lemma 1 shows that a submanifold is Spin$^c$ if and only if its normal bundle is. In [2] it was observed that therefore all Lagrangian submanifolds in Kähler manifolds of complex dimension $\leq 4$ have Spin$^c$ normal bundles. The following Proposition shows that this is sharp.

**Proposition 1.** *For every $n \geq 5$ there is a compact Kähler manifold $M$ of complex dimension $n$ with a Lagrangian submanifold $Z$ whose normal bundle is not* Spin$^c$.

*Proof.* It suffices to prove the case $n = 5$, as we can then take products with Abelian varieties (and Lagrangian subtori) to get the general case.

Consider the non-Spin$^c$ 5-manifold $Z = \mathrm{SU}(3)/\mathrm{SO}(3)$, cf. the proof of Lemma 7. This is the real part of its complexification $A$, which is an affine algebraic variety homogeneous under the complexification of $\mathrm{SU}(3)$. Let $M$ be a smooth projective variety obtained as (a resolution of) a projective closure of $A$. Then $Z$ is Lagrangian with respect to the Fubini-Study metric of $M$, and does not have Spin$^c$ normal bundle by Lemma 1. □

We do not know whether such examples exist in Calabi-Yau manifolds.

Mathematisches Institut, Ludwig-Maximilians-Universität München, Theresienstr. 39, 80333 München, Germany

*E-mail address*: `bohr@rz.mathematik.uni-muenchen.de`

*E-mail address*: `hanke@rz.mathematik.uni-muenchen.de`

*E-mail address*: `dieter@member.ams.org`